# Poisson-type deviation inequalities for curved continuous-time Markov chains

ALDÉRIC JOULIN

*MODAL'X, Bât. G, Université Paris 10, 200, Avenue de la République, 92001 Nanterre Cedex, France. E-mail: ajoulin@u-paris10.fr*

In this paper, we present new Poisson-type deviation inequalities for continuous-time Markov chains whose Wasserstein curvature or $\Gamma$-curvature is bounded below. Although these two curvatures are equivalent for Brownian motion on Riemannian manifolds, they are not comparable in discrete settings and yield different deviation bounds. In the case of birth–death processes, we provide some conditions on the transition rates of the associated generator for such curvatures to be bounded below and we extend the deviation inequalities established [Ané, C. and Ledoux, M. On logarithmic Sobolev inequalities for continuous time random walks on graphs. *Probab. Theory Related Fields* **116** (2000) 573–602] for continuous-time random walks, seen as models in null curvature. Some applications of these tail estimates are given for Brownian-driven Ornstein–Uhlenbeck processes and $M/M/1$ queues.

## 1. Introduction

Let $\mu$ be a probability measure on a metric space $(E, d)$ and let $h: \mathbb{R}_+ \to \mathbb{R}_+$ be a function tending to 0 at infinity. The measure $\mu$ is said to satisfy a deviation inequality of speed $h$ if, for any real Lipschitz function $f$ on $(E, d)$ with Lipschitz constant smaller than 1, the following inequality holds:

$$\mu(f - \mu(f) \geq y) \leq h(y), \qquad y > 0,$$

where $\mu(f) = \int f \, d\mu$. Applying the above inequality to $-f$ also leads to a concentration inequality stating that any Lipschitz map is concentrated around its mean under $\mu$ with high probability. In particular, the concentration is Gaussian if $h$ is of order $\exp(-y^2)$, whereas it is of Poisson type if $h$ is of order $\exp(-y \log(y))$ for large $y$.

Actually, the concentration-of-measure phenomenon is useful for determining the rate of convergence of functionals involving a large number of random variables. In recent years, this area has been extensively investigated in the context of dependent random variables such as Markov chains. For instance, Gaussian concentration was considered by Marton [11] and Djellout *et al.* [6] by means of transportation cost inequalities, whereas Ané and Ledoux [1], Samson [12] and Houdré and Tetali [9] established some appropriate functional inequalities to derive Gaussian and Poisson-type deviation inequalities.

The purpose of the present paper is to give new deviation inequalities of Poisson type for continuous-time Markov chains, which extend the tail estimates of Ané and Ledoux [1]. Our approach is based on semigroup analysis and uses the notion of curvature for Markov processes on general metric measure spaces recently investigated by Sturm and Von Renesse [14]. Although the various Brownian curvatures on a smooth Riemannian manifold are essentially equivalent and characterize the uniform lower bounds on the Ricci curvature of the manifold, such an equivalence does not hold for continuous-time Markov chains since in general, discrete gradients do not satisfy the chain rule formula. Thus it is natural to study the role played by each type of curvature in the concentration-of-measure phenomenon.

The paper is organized as follows. In Section 2, two different notions of curvature of continuous-time Markov chains are introduced: Wasserstein curvature and $\Gamma$-curvature. Section 3 is concerned with the main results of the paper. Namely, a Poisson-type deviation inequality is established in Theorem 3.1 for continuous-time Markov chains with bounded angle bracket, provided the Wasserstein curvature of the process is bounded below. Under the hypothesis of a lower bound on the $\Gamma$-curvature, a general estimate is derived in Theorem 3.4 and with further assumptions on the chain, the latter upper bound is computed to yield Poisson tail probabilities for processes with non-necessarily bounded angle bracket. The case of birth–death processes on $\mathbb{N}$ or $\{0, 1, \ldots, n\}$ is investigated in Section 4, where we give some conditions on the transition rates of the associated generator for such curvatures to be bounded below. As a corollary of the tail estimates emphasized above, we extend to birth–death processes the deviation inequalities established by Ané and Ledoux [1] for continuous-time random walks on graphs, seen as models in null curvature. Finally, some applications of these tail estimates are given to Brownian-driven Ornstein–Uhlenbeck processes and $M/M/1$ queues.

## 2. Notation and preliminaries on curvatures

Throughout this paper, $E$ is a countable set endowed with a metric $d$ which differs from the trivial one $\varrho(x,y) = 1_{x \neq y}$, $x, y \in E$, $\mathscr{F}(E)$ is the collection of all real-valued functions on $E$, $\mathscr{B}(E) \subset \mathscr{F}(E)$ is the subspace of bounded functions on $E$ equipped with the supremum norm $\|f\|_\infty = \sup_{x \in E} |f(x)|$ and the space $\mathrm{Lip}_d(E)$ consists of Lipschitz function $f$ on $E$ with Lipschitz seminorm

$$\|f\|_{\mathrm{Lip}_d} := \sup_{x \neq y} \frac{|f(x) - f(y)|}{d(x,y)} < +\infty.$$

On a filtered probability space $(\Omega, \mathscr{F}, (\mathscr{F}_t)_{t \geq 0}, \mathbb{P})$, let $(X_t)_{t \geq 0}$ be an $E$-valued regular non-explosive continuous-time Markov chain, where regularity is understood in the sense of Chen [5]. The generator $\mathcal{L}$ of the chain is given for any $f \in \mathscr{B}(E)$ by

$$\mathcal{L}f(x) = \sum_{y \in E}(f(y) - f(x))Q(x,y), \qquad x \in E,$$

where the transition rates $(Q(x,y))_{x \neq y}$ are non-negative. Let $(P_t)_{t \geq 0}$ be the associated semigroup which acts on the elements of $\mathscr{B}(E)$ as follows:

$$P_t f(x) := \mathbb{E}_x[f(X_t)] = \sum_{y \in E} f(y) P_t(x,y), \qquad x \in E.$$

Denote by $\mathscr{P}_1(E)$ the space of probability measures $\mu$ on $E$ such that $\sum_{y \in E} d(y,z)\mu(y) < +\infty$ for some $z \in E$. If the Markov kernel $P_t(x,\cdot) \in \mathscr{P}_1(E)$ for some $x \in E$, then the semigroup $(P_t)_{t \geq 0}$ is also well defined on the space $\mathrm{Lip}_d(E)$.

If there exists a constant $V > 0$ such that $\|\sum_{y \in E} d(\cdot,y)^2 Q(\cdot,y)\|_\infty \leq V^2$, we say that the angle bracket of the process $(X_t)_{t \geq 0}$ is bounded by $V^2$. Moreover, the jumps of $(X_t)_{t \geq 0}$ are said to be bounded by a positive constant $b > 0$ if $\sup_{t > 0} d(X_{t-}, X_t) \leq b$.

Let us introduce the notion of curved Markov chains in the Wasserstein sense.

**Definition 2.1.** *Assume that the Markov kernel $P_t(x,\cdot) \in \mathscr{P}_1(E)$ for some $x \in E$. The $d$-Wasserstein curvature at time $t > 0$ of the continuous-time Markov chain $(X_t)_{t \geq 0}$ is defined by*

$$K_t := -\frac{1}{t} \sup\left\{\log\left(\frac{\|P_t f\|_{\mathrm{Lip}_d}}{\|f\|_{\mathrm{Lip}_d}}\right) : f \in \mathrm{Lip}_d(E), f \neq \mathrm{const}\right\} \in [-\infty, +\infty).$$

*It is said to be bounded below by $K \in \mathbb{R}$ if $\inf_{t > 0} K_t \geq K$. In other words, for any Lipschitz function $f \in \mathrm{Lip}_d(E)$ and any $t > 0$,*

$$\|P_t f\|_{\mathrm{Lip}_d} \leq e^{-Kt} \|f\|_{\mathrm{Lip}_d}. \qquad (2.1)$$

**Remark 2.2.** Given $\mu, \nu \in \mathscr{P}_1(E)$, define the $d$-Wasserstein distance between $\mu$ and $\nu$ by

$$W_d(\mu, \nu) := \inf_\pi \sum_{x,y \in E} d(x,y)\pi(x,y),$$

where the infimum runs over all $\pi \in \mathscr{P}_1(E \times E)$ with marginals $\mu$ and $\nu$. By the Kantorovich–Rubinstein duality theorem (cf. Chen [5], Theorem 5.10), the $d$-Wasserstein distance can be rewritten as

$$W_d(\mu, \nu) = \sup\left\{\left|\sum_{x \in E} f(x)(\mu(x) - \nu(x))\right| : \|f\|_{\mathrm{Lip}_d} \leq 1\right\}.$$

Hence, the following assertions are equivalent:

(i) $\inf_{t > 0} K_t \geq K$;

(ii) $W_d(P_t(x,\cdot), P_t(y,\cdot)) \leq e^{-Kt} d(x,y)$, for any $x, y \in E$ and any $t > 0$.

Note that assertion (ii) characterizes the lower bounds on the $d$-Wasserstein curvature in terms of contraction properties of the semigroup in the metric $W_d$, which induces a coupling approach. Such an inequality was introduced by Marton [11], with the trivial metric $\varrho$, and also by Djellout *et al.* [6] through condition (C1), to establish transportation and Gaussian concentration inequalities for weakly dependent sequences.

In order to introduce the "carré du champ" operator, we follow Bakry [2] and assume the existence of an algebra, say $\mathcal{A}$, which contains the bounded functions and is stable under the action of $\mathcal{L}$, $P_t$ and by composition with the $C^\infty$-functions. The "carré du champ" operator $\Gamma$ is defined on $\mathcal{A} \times \mathcal{A}$ by

$$\Gamma(f,g)(x) := \tfrac{1}{2}(\mathcal{L}(fg)(x) - f(x)\mathcal{L}g(x) - g(x)\mathcal{L}f(x))$$
$$= \tfrac{1}{2} \sum_{y \in E} (f(y) - f(x))(g(y) - g(x))Q(x,y).$$

We set $\Gamma f = \Gamma(f,f)$ and introduce the notion of curved Markov chains in the $\Gamma$-sense.

**Definition 2.3.** *The $\Gamma$-curvature at time $t > 0$ of the continuous-time Markov chain $(X_t)_{t \geq 0}$ is defined by*

$$\rho_t := -\frac{1}{t} \sup \left\{ \log \left\| \frac{(\Gamma P_t f)^{1/2}}{P_t(\Gamma f)^{1/2}} \right\|_\infty : f \in \mathcal{A}, f \neq \mathrm{const} \right\} \in [-\infty, +\infty).$$

*It is said to be bounded below by $\rho \in \mathbb{R}$ if $\inf_{t>0} \rho_t \geq \rho$. In other words, for any $f \in \mathcal{A}$, any $x \in E$ and any $t > 0$,*

$$(\Gamma P_t f)^{1/2}(x) \leq e^{-\rho t} P_t(\Gamma f)^{1/2}(x). \tag{2.2}$$

**Remark 2.4.** Such an inequality is the discrete analogue of the commutation relation between local gradient and heat kernel on Riemannian manifolds with Ricci curvature bounded below; see Bakry and Émery [3].

As mentioned in the Introduction, the two curvatures are equivalent in the continuous setting of Brownian motions on Riemannian manifolds (cf. Sturm and Von Renesse [14]). This equivalence does not hold for discrete spaces since in general, the discrete gradients do not satisfy the chain rule formula and the curvatures defined above are not comparable. However, we point out that the semigroup appears in both sides of the inequality (2.2), whereas it is dropped in the right-hand side of inequality (2.1). Hence, we deduce that the assumption (2.2) is stronger than (2.1) in some sense, and we expect to obtain stronger deviation results when dealing with $\Gamma$-curvature.

To complete the preliminaries, let us make some comments on the deviation inequalities we will establish in the remainder of this paper.

1) Our estimates are given for the distribution of $X_t$ given $X_0 = x$, uniformly in $x \in E$ and for any $t > 0$. Hence, without risk of confusion, the ranges of validity of the parameters $x$ and $t$ will not be mentioned in our results.
2) In order to ease notation, our results are given with the function $u \mapsto u \log(1+u)/2$ in the upper bounds. However, sharper estimates are also available when replacing this function by $u \mapsto (1+u)\log(1+u) - u$, $u \geq 0$.

## 3. Deviation bounds for curved continuous-time Markov chains

In this section, we present Poisson-type deviation estimates under the assumption of a lower bound on the curvatures of the continuous-time Markov chain $(X_t)_{t \geq 0}$. Let us start with $d$-Wasserstein curvature.

**Theorem 3.1.** *Assume that the jumps and the angle bracket of $(X_t)_{t \geq 0}$ are bounded by $b > 0$ and $V^2 > 0$, respectively. Suppose, moreover, that its $d$-Wasserstein curvature is bounded below by $K \in \mathbb{R}$. Let $f \in \mathrm{Lip}_d(E)$ and define $C_{t,K} := \sup_{0 \leq s \leq t} e^{-K(t-s)}$ and $M_{t,K} := (1 - e^{-2Kt})/(2K)$ ($M_{t,K} = t$ if $K = 0$). Then, for any $y > 0$,*

$$\mathbb{P}_x(f(X_t) - \mathbb{E}_x[f(X_t)] \geq y)$$
$$\leq \exp\left(-\frac{y}{2bC_{t,K}\|f\|_{\mathrm{Lip}_d}} \log\left(1 + \frac{bC_{t,K}y}{M_{t,K}V^2\|f\|_{\mathrm{Lip}_d}}\right)\right). \tag{3.1}$$

**Proof.** First assume that the Lipschitz function $f$ is bounded. The process $(Z_s^f)_{0 \leq s \leq t}$ given by $Z_s^f := P_{t-s}f(X_s) - P_t f(X_0)$ is then a real $\mathbb{P}_x$-martingale with respect to the truncated filtration $(\mathscr{F}_s)_{0 \leq s \leq t}$ and we have, by Itô's formula,

$$Z_s^f = \sum_{y,z \in E} \int_0^s (P_{t-\tau}f(y) - P_{t-\tau}f(z))1_{\{X_{\tau-}=z\}}(N_{z,y} - \sigma_{z,y})(\mathrm{d}\tau),$$

where $(N_{z,y})_{z,y \in E}$ is a family of independent Poisson processes on $\mathbb{R}_+$ with respective intensity $\sigma_{z,y}(\mathrm{d}t) = Q(z,y)\,\mathrm{d}t$. Since the $d$-Wasserstein curvature is bounded below, the jumps of $(Z_s^f)_{0 \leq s \leq t}$ are bounded:

$$|Z_s^f - Z_{s-}^f| = |P_{t-s}f(X_s) - P_{t-s}f(X_{s-})|$$
$$\leq d(X_s, X_{s-})\|f\|_{\mathrm{Lip}_d} C_{t,K}$$
$$\leq b\|f\|_{\mathrm{Lip}_d} C_{t,K}.$$

Moreover, the angle bracket is also bounded:

$$\langle Z^f, Z^f \rangle_s = \sum_{y,z \in E} \int_0^s (P_{t-\tau}f(y) - P_{t-\tau}f(z))^2 1_{\{X_{\tau-}=z\}} \sigma_{z,y}(\mathrm{d}\tau)$$

$$\leq \|f\|_{\mathrm{Lip}_d}^2 \sum_{y,z \in E} \int_0^s \mathrm{e}^{-2K(t-\tau)} d(z,y)^2 \mathbf{1}_{\{X_{\tau^-}=z\}} Q(z,y) \, \mathrm{d}\tau$$

$$\leq \|f\|_{\mathrm{Lip}_d}^2 M_{t,K} V^2.$$

By Kallenberg [10], Lemma 23.19, for any positive $\lambda$, the process $(Y_s^{(\lambda)})_{0 \leq s \leq t}$ given by

$$Y_s^{(\lambda)} := \exp\{\lambda Z_s^f - \lambda^2 \psi(\lambda b \|f\|_{\mathrm{Lip}_d} C_{t,K}) \langle Z^f, Z^f \rangle_s\}$$

is a $\mathbb{P}_x$-supermartingale with respect to $(\mathscr{F}_s)_{0 \leq s \leq t}$, where $\psi(z) = z^{-2}(\mathrm{e}^z - z - 1)$, $z > 0$. Thus we get, for any $\lambda > 0$,

$$\mathbb{E}_x[\exp\{\lambda(f(X_t) - \mathbb{E}_x[f(X_t)])\}]$$
$$= \mathbb{E}_x[\mathrm{e}^{\lambda Z_t^f}]$$
$$\leq \exp\{\lambda^2 \|f\|_{\mathrm{Lip}_d}^2 M_{t,K} V^2 \psi(\lambda b \|f\|_{\mathrm{Lip}_d} C_{t,K})\} \mathbb{E}_x[Y_t^{(\lambda)}]$$
$$\leq \exp\{\lambda^2 \|f\|_{\mathrm{Lip}_d}^2 M_{t,K} V^2 \psi(\lambda b \|f\|_{\mathrm{Lip}_d} C_{t,K})\}$$
$$= \exp\left\{\frac{M_{t,K} V^2}{b^2 C_{t,K}^2} (\mathrm{e}^{\lambda b \|f\|_{\mathrm{Lip}_d} C_{t,K}} - \lambda b \|f\|_{\mathrm{Lip}_d} C_{t,K} - 1)\right\}.$$

Using Chebyshev's inequality and optimizing in $\lambda > 0$ in the exponential estimate above, the deviation inequality (3.1) is then established in the bounded case. Finally, the boundedness assumption on $f$ is removed by a classical argument. □

**Remark 3.2.** If $K = 0$, then the estimate in Theorem 3.1 is similar to the deviation inequalities established by Houdré [8] and Schmuckenschläger [13] for infinitely divisible distributions with compactly supported Lévy measure. If $K < 0$, the decay in (3.1) is slower, due to some exponential factors, whereas if $K > 0$, the chain is ergodic (cf. Chen [5], Theorem 5.23), and such an estimate can be extended, as $t$ goes to infinity, to the stationary distribution, as illustrated in Section 4. On the other hand, the sign of $K$ has no influence in small time on (3.1).

Note that Theorem 3.1 does not allow us to consider continuous-time Markov chains with unbounded angle bracket. To overcome this difficulty, we must impose some assumptions on a different curvature of the process, namely $\Gamma$-curvature.

Our present purpose is to adapt to the Markovian case the covariance method of Houdré [8] in order to derive deviation inequalities for curved continuous-time Markov chains in the $\Gamma$-sense. Although Wasserstein curvature and $\Gamma$-curvature are not comparable in discrete spaces, the results we give now are more general than in Theorem 3.1.

Before turning to Theorem 3.4 below, let us establish the following.

**Lemma 3.3.** *Assume that $(X_t)_{t \geq 0}$ has $\Gamma$-curvature bounded below by $\rho \in \mathbb{R}$. Let $g_1, g_2 \in \mathscr{B}(E)$ with $\|\Gamma g_1\|_\infty < +\infty$ and define $L_{t,\rho} = (1 - \mathrm{e}^{-2\rho t})/(2\rho)$ if $\rho \neq 0$ and $L_{t,\rho} = t$ other-* wise. We then have the covariance inequality

$$\mathrm{Cov}_x[g_1(X_t), g_2(X_t)] \leq 2L_{t,\rho}\|\Gamma g_1\|_\infty^{1/2}\mathbb{E}_x[(\Gamma g_2)^{1/2}(X_t)].$$

**Proof.** As in the proof of Theorem 3.1, we have, for $i = 1, 2$,

$$g_i(X_t) - \mathbb{E}_x[g_i(X_t)] = \sum_{y,z \in E} \int_0^t (P_{t-s}g_i(y) - P_{t-s}g_i(z))1_{\{X_{s-}=z\}}(N_{z,y} - \sigma_{z,y})(\mathrm{d}s).$$

By the Cauchy–Schwarz inequality,

$$\begin{aligned}
\mathrm{Cov}_x[g_1(X_t), g_2(X_t)] &= 2\int_0^t P_s(\Gamma(P_{t-s}g_1, P_{t-s}g_2))(x)\,\mathrm{d}s \\
&\leq 2\int_0^t P_s((\Gamma P_{t-s}g_1)^{1/2}(\Gamma P_{t-s}g_2)^{1/2})(x)\,\mathrm{d}s \\
&\leq 2\int_0^t \exp\{-2\rho(t-s)\}P_s(P_{t-s}(\Gamma g_1)^{1/2}P_{t-s}(\Gamma g_2)^{1/2})(x)\,\mathrm{d}s,
\end{aligned}$$

where, in the latter inequality, we used the assumption of a lower bound $\rho$ on the $\Gamma$-curvature. Since $(P_t)_{t\geq 0}$ is a contraction operator on $\mathscr{B}(E)$, we have

$$\begin{aligned}
\mathrm{Cov}_x[g_1(X_t), g_2(X_t)] &\leq 2\|\Gamma g_1\|_\infty^{1/2}\int_0^t \exp\{-2\rho(t-s)\}P_s(P_{t-s}(\Gamma g_2)^{1/2})(x)\,\mathrm{d}s \\
&= 2L_{t,\rho}\|\Gamma g_1\|_\infty^{1/2}\mathbb{E}_x[(\Gamma g_2)^{1/2}(X_t)].
\end{aligned}$$

The proof is complete. $\square$

We are now able to state Theorem 3.4, which presents a general deviation bound for curved continuous-time Markov chains in the $\Gamma$-sense.

**Theorem 3.4.** *Assume that $(X_t)_{t\geq 0}$ has $\Gamma$-curvature bounded below by $\rho \in \mathbb{R}$. Let $f \in \mathcal{A} \cap \mathrm{Lip}_d(E)$ with $\|\Gamma f\|_\infty < +\infty$ and define the function $\psi_{f,t}:\mathbb{R}_+ \to \mathbb{R}_+ \cup \{\infty\}$ by*

$$\psi_{f,t}(\lambda) := \sqrt{2}L_{t,\rho}\|\Gamma f\|_\infty^{1/2}\left\|\sum_{y\in E}(f(y) - f(\cdot))^2\left(\frac{\exp\{\lambda\|f\|_{\mathrm{Lip}_d}d(\cdot,y)\} - 1}{\|f\|_{\mathrm{Lip}_d}d(\cdot,y)}\right)^2 Q(\cdot,y)\right\|_\infty^{1/2},$$

*where $L_{t,\rho}$ is defined in Lemma 3.3. Letting $M_{f,t} := \sup\{\lambda > 0 : \psi_{f,t}(\lambda) < +\infty\}$, we have*

$$\mathbb{P}_x(f(X_t) - \mathbb{E}_x[f(X_t)] \geq y) \leq \exp\inf_{\lambda \in (0, M_{f,t})}\int_0^\lambda (\psi_{f,t}(\tau) - y)\,\mathrm{d}\tau, \qquad y > 0. \quad (3.2)$$

**Remark 3.5.** Note that $\psi_{f,t}$ is bijective from $(0, M_{f,t})$ to $(0, +\infty)$, so the term in the exponential is negative and inequality (3.2) makes sense.

**Proof.** Using a standard argument, we can reduce the problem to establishing the result for $f$ bounded Lipschitz. Applying the covariance inequality of Lemma 3.3 with the functions $g_1(z) = f(z) - \mathbb{E}_x[f(X_t)]$ and $g_2(z) = \exp\{\lambda(f(z) - \mathbb{E}_x[f(X_t)])\}$, $z \in E$, $\lambda \in (0, M_{f,t})$, we have

$$\mathbb{E}_x[(f(X_t) - \mathbb{E}_x[f(X_t)]) \exp\{\lambda(f(X_t) - \mathbb{E}_x[f(X_t)])\}]$$
$$\leq 2L_{t,\rho} \|\Gamma f\|_\infty^{1/2} \exp\{-\lambda \mathbb{E}_x[f(X_t)]\} \mathbb{E}_x[(\Gamma e^{\lambda f})^{1/2}(X_t)]$$
$$\leq \sqrt{2} L_{t,\rho} \|\Gamma f\|_\infty^{1/2} \mathbb{E}_x \left[ \exp\{\lambda(f(X_t) - \mathbb{E}_x[f(X_t)])\} \right.$$
$$\left. \times \left( \sum_{y,z \in E} (\exp\{\lambda|f(y) - f(z)|\} - 1)^2 1_{\{X_t=z\}} Q(z,y) \right)^{1/2} \right]$$
$$\leq \psi_{f,t}(\lambda) \mathbb{E}_x[\exp\{\lambda(f(X_t) - \mathbb{E}_x[f(X_t)])\}].$$

Letting $H_{f,t,x}(\lambda) := \mathbb{E}_x[\exp\{\lambda(f(X_t) - \mathbb{E}_x[f(X_t)])\}]$, the latter inequality can be rewritten as $H'_{f,t,x}(\lambda) \leq \psi_{f,t}(\lambda) H_{f,t,x}(\lambda)$, from which follows the bound

$$\mathbb{E}_x[\exp\{\lambda(f(X_t) - \mathbb{E}_x[f(X_t)])\}] = H_{f,t,x}(\lambda) \leq \exp\left\{\int_0^\lambda \psi_{f,t}(\tau) \, d\tau\right\}, \qquad \lambda \in (0, M_{f,t}).$$

Finally, using Chebyshev's inequality, Theorem 3.4 is established. □

Since the estimate (3.2) is very general, let us make further assumptions on the process $(X_t)_{t\geq 0}$ in order to get Poisson-type deviation inequalities. Denote in the sequel $L_{t,\rho} = (1 - e^{-2\rho t})/(2\rho)$ if $\rho \neq 0$ and $L_{t,\rho} = t$ otherwise. Using the notation of Theorem 3.4, we have the following.

**Corollary 3.6.** *Under the hypothesis of Theorem 3.4, suppose further that the jumps of $(X_t)_{t\geq 0}$ are bounded by $b > 0$. Then, for any $y > 0$,*

$$\mathbb{P}_x(f(X_t) - \mathbb{E}_x[f(X_t)] \geq y) \leq \exp\left(-\frac{y}{2b\|f\|_{\mathrm{Lip}_d}} \log\left(1 + \frac{yb\|f\|_{\mathrm{Lip}_d}}{2L_{t,\rho}\|\Gamma f\|_\infty}\right)\right).$$

**Proof.** Under the notation of Theorem 3.4, the boundedness of the jumps implies $M_{f,t} = +\infty$, and $\psi_{f,t}$ is bounded by

$$\psi_{f,t}(\lambda) \leq 2L_{t,\rho} \|\Gamma f\|_\infty \frac{\exp\{\lambda b \|f\|_{\mathrm{Lip}_d}\} - 1}{b\|f\|_{\mathrm{Lip}_d}}, \qquad \lambda > 0.$$

Using Theorem 3.4 and optimizing in $\lambda > 0$, the proof is then achieved. □

Note that the latter deviation inequality is more general than (3.1) since the finiteness assumption on $\|\Gamma f\|_\infty$ allows us to relax the boundedness assumption on the angle

bracket of the process $(X_t)_{t\geq 0}$. Thus when the angle bracket is bounded, the next corollary exhibits an estimate comparable to that of Theorem 3.1.

**Corollary 3.7.** *Assume that the jumps and the angle bracket of $(X_t)_{t\geq 0}$ are bounded by $b > 0$ and $V^2 > 0$, respectively. Suppose further that its $\Gamma$-curvature is bounded below by $\rho \in \mathbb{R}$. Letting $f \in \mathcal{A} \cap \mathrm{Lip}_d(E)$, then for any $y > 0$,*

$$\mathbb{P}_x(f(X_t) - \mathbb{E}_x[f(X_t)] \geq y) \leq \exp\left(-\frac{y}{2b\|f\|_{\mathrm{Lip}_d}}\log\left(1 + \frac{by}{L_{t,\rho}V^2\|f\|_{\mathrm{Lip}_d}}\right)\right).$$

**Proof.** By the boundedness of the jumps and of the angle bracket, the function $\psi_{f,t}$ in Theorem 3.4 is bounded by

$$\psi_{f,t}(\lambda) \leq L_{t,\rho}V^2\|f\|_{\mathrm{Lip}_d}\frac{\exp\{\lambda b\|f\|_{\mathrm{Lip}_d}\} - 1}{b}, \qquad \lambda > 0.$$

Finally, applying Theorem 3.4 yields the result. □

## 4. The case of birth–death processes

In the paper of Ané and Ledoux [1], some deviation inequalities are established for continuous-time random walks on graphs. Such processes may be seen as models in null curvature since the transition rates of the associated generator do not depend on the space variable. Using the results of Section 3, the purpose of this section is to extend these tail estimates to birth–death processes whose curvatures are bounded below.

Let $(X_t)_{t\geq 0}$ be a regular non-explosive birth–death process with stationary distribution $\pi$ on the state space $E = \mathbb{N}$ or $E = \{0, 1, \ldots, n\}$, endowed with the classical metric $d(x, y) = |x - y|$, $x, y \in E$. Such a process is a continuous-time Markov chain with generator defined on $\mathscr{F}(E)$ by

$$\mathcal{L}f(x) = \lambda_x(f(x+1) - f(x)) + \nu_x(f(x-1) - f(x)), \qquad x \in E, \qquad (4.1)$$

where the transition rates $\lambda$ and $\nu$ are positive with 0 as reflecting state, that is, $\nu_0 = 0$ (if $E = \{0, 1, \ldots, n\}$, the state $n$ is also reflecting: $\lambda_n = 0$), ensuring irreducibility. Denote by $(P_t)_{t\geq 0}$ the homogeneous semigroup whose transition probabilities are given for any $x \in E$ by

$$P_t(x, y) = \begin{cases} \lambda_x t + o(t), & \text{if } y = x + 1, \\ \nu_x t + o(t), & \text{if } y = x - 1, \\ 1 - (\lambda_x + \nu_x)t + o(t), & \text{if } y = x, \end{cases}$$

where the function $o$ is such that $o(t)/t$ converges to 0 as $t$ tends to 0. We assume in the remainder of the paper that the stationary distribution $\pi \in \mathscr{P}_1(E)$. Then, by the contraction property of $P_t$ in $L^1(\pi)$, the Markov kernel $P_t(x, \cdot)$ belongs to the space $\mathscr{P}_1(E)$ for any $x \in E$ and the semigroup is well defined on the space $\mathrm{Lip}_d(E)$.

In the case of birth–death processes, the "carré du champ" operator $\Gamma$ is given on $\mathcal{A} = \mathscr{F}(E)$ by

$$\Gamma f(x) = \tfrac{1}{2}\{\lambda_x(f(x+1) - f(x))^2 + \nu_x(f(x-1) - f(x))^2\}, \qquad x \in E.$$

### 4.1. Criteria for lower bounded curvatures

Let us give some criteria on the generator of the process $(X_t)_{t \geq 0}$ which ensure that the different curvatures are bounded below.

**Proposition 4.1.** *Assume that there exists a real number $K$ such that the transition rates $\lambda$ and $\nu$ satisfy the inequality*

$$\inf_{x \in E \setminus \{0\}} \lambda_{x-1} - \lambda_x + \nu_x - \nu_{x-1} \geq K. \qquad (4.2)$$

*Then the $d$-Wasserstein curvature of the process $(X_t)_{t \geq 0}$ is bounded below by $K$.*

**Proof.** Let us establish the result via a coupling argument. Consider $(X_t^x)_{t \geq 0}$ and $(X_t^y)_{t \geq 0}$, two independent copies of the process $(X_t)_{t \geq 0}$, starting from $x$ and $y$, respectively. The generator $\tilde{\mathcal{L}}$ of the process $(X_t^x, X_t^y)_{t \geq 0}$ is then given, for any $f \in \mathscr{F}(E \times E)$, by

$$\tilde{\mathcal{L}} f(z, w) = (\mathcal{L} f(\cdot, w))(z) + (\mathcal{L} f(z, \cdot))(w), \qquad z, w \in E.$$

Since the transition rates of the generator satisfy (4.2), we immediately have the bound $\tilde{\mathcal{L}} d(z, z+1) \leq -K$, $z \in E$, which is equivalent to the inequality $\tilde{\mathcal{L}} d(z, w) \leq -K d(z, w)$ for any $z, w \in E$. Therefore, we obtain the estimate $\mathbb{E}[d(X_t^x, X_t^y)] \leq e^{-Kt} d(x, y)$, which, in turn, implies

$$W_d(P_t(x, \cdot), P_t(y, \cdot)) \leq e^{-Kt} d(x, y).$$

Finally, by the equivalent statements of Remark 2.2, the $d$-Wasserstein curvature of $(X_t)_{t \geq 0}$ is bounded below by $K$. $\square$

**Remark 4.2.** If $E = \mathbb{N}$ and the transition rates of the generator are bounded and satisfy (4.2), then necessarily $K \leq 0$.

In order to establish modified logarithmic Sobolev inequalities for continuous-time random walks on $\mathbb{Z}$, Ané and Ledoux [1] used a suitable $\Gamma_2$-calculus to give a criterion under which the $\Gamma$-curvature is bounded below by 0. Actually, this criterion can be generalized to any real lower bound on the $\Gamma$-curvature via Lemma 4.3 below.

Define the operator $\Gamma_2$ on $\mathscr{F}(E)$ by

$$\Gamma_2 f(x) := \tfrac{1}{2}(\mathcal{L} \Gamma f(x) - 2\Gamma(f, \mathcal{L} f)(x)), \qquad x \in E.$$

By adapting the proof of Ané and Ledoux [1] mentioned above, we obtain the following lemma.

**Lemma 4.3.** *Assume that there exists $\rho \in \mathbb{R}$ such that for any $f \in \mathscr{F}(E)$,*

$$\Gamma_2 f(x) - \Gamma(\Gamma f)^{1/2}(x) \geq \rho \Gamma f(x), \qquad x \in E. \tag{4.3}$$

*Then the $\Gamma$-curvature of the process $(X_t)_{t \geq 0}$ is bounded below by $\rho$.*

**Proposition 4.4.** *Assume that there exists some $\rho \geq 0$ such that the transition rates $\lambda$ and $\nu$ satisfy*

$$\inf_{x \in E \setminus \{0, \sup E\}} \min\{\lambda_{x-1} - \lambda_x, \nu_{x+1} - \nu_x\} \geq \rho. \tag{4.4}$$

*Then the $\Gamma$-curvature of the process $(X_t)_{t \geq 0}$ is bounded below by $\rho$.*

**Proof.** By Lemma 4.3, the result holds true if the inequality (4.3) above is satisfied. Let us prove this inequality.

Fix $x \in E$ and let $a = f(x) - f(x+1)$, $b = f(x) - f(x-1)$, $c = f(x+2) - f(x+1)$ and $d = f(x-2) - f(x-1)$. We have

$$2\Gamma_2 f(x) - 2\Gamma(\Gamma f)^{1/2}(x) = \lambda_x(\nu_{x+1} - \nu_x)a^2 + \nu_x(\lambda_{x-1} - \lambda_x)b^2 + I(x) + J(x),$$

where

$$I(x) := \lambda_x \lambda_{x+1} ac + \lambda_x \nu_x ab + \lambda_x (\lambda_{x+1} c^2 + \nu_{x+1} a^2)^{1/2} (\lambda_x a^2 + \nu_x b^2)^{1/2},$$

$$J(x) := \nu_x \nu_{x-1} bd + \lambda_x \nu_x ab + \nu_x (\lambda_{x-1} b^2 + \nu_{x-1} d^2)^{1/2} (\lambda_x a^2 + \nu_x b^2)^{1/2}.$$

Since the transition rates $\lambda$ and $\nu$ satisfy (4.4), we get

$$2\Gamma_2 f(x) - 2\Gamma(\Gamma f)^{1/2}(x) \geq 2\rho \Gamma f(x) + I(x) + J(x).$$

By symmetry with the function $J$, it is sufficient to establish $I \geq 0$. We have

$$I(x) \geq \lambda_x (\lambda_{x+1} c^2 + \nu_{x+1} a^2)^{1/2} (\lambda_x a^2 + \nu_x b^2)^{1/2} - \lambda_x \lambda_{x+1} |ac| - \lambda_x \nu_x |ab|$$
$$= \lambda_x (I_1(x) - I_2(x)),$$

where

$$I_1(x) := (\lambda_{x+1} c^2 + \nu_{x+1} a^2)^{1/2} (\lambda_x a^2 + \nu_x b^2)^{1/2} \quad \text{and} \quad I_2(x) := \lambda_{x+1} |ac| + \nu_x |ab|.$$

Again using inequality (4.4),

$$(I_1(x))^2 - (I_2(x))^2$$
$$= \lambda_{x+1}(\lambda_x - \lambda_{x+1})a^2 c^2 + \nu_x(\nu_{x+1} - \nu_x)a^2 b^2 + \lambda_x \nu_{x+1} a^4 + \lambda_{x+1} \nu_x b^2 c^2 - 2\nu_x \lambda_{x+1} a^2 bc$$
$$\geq \nu_x \lambda_{x+1} (a^2 - bc)^2 \geq 0.$$

The proof is complete. $\square$

*Remark 4.5.* We note that the equivalence holds in Lemma 4.3. Indeed, if the $\Gamma$-curvature of the process $(X_t)_{t\geq 0}$ is bounded below by $\rho$, then for any $f \in \mathscr{F}(E)$, the function $\alpha$ defined on $[0, \infty)$ by $\alpha(t) = \mathrm{e}^{-\rho t} P_t \sqrt{\Gamma f} - \sqrt{\Gamma P_t f}$ is non-negative with $\alpha(0) = 0$. Hence, we have $\alpha'(0) \geq 0$, which is the inequality (4.3).

Moreover, we point out that if $E = \mathbb{N}$ and the transition rates of the generator satisfy (4.4), then necessarily $\rho = 0$.

## 4.2. Applications

The proofs of the following results are omitted since they are immediate applications of Theorem 3.1 (resp., Theorem 3.4), once the assumptions of Proposition 4.1 (resp., Proposition 4.4) are satisfied. Let us start with the case $E = \mathbb{N}$.

**Corollary 4.6.** *Assume that the transition rates $\lambda$ and $\nu$ are bounded on $\mathbb{N}$ and suppose that there exists $K \leq 0$ such that $\inf_{x \in \mathbb{N} \setminus \{0\}} \lambda_{x-1} - \lambda_x + \nu_x - \nu_{x-1} \geq K$. If $f \in \mathrm{Lip}_d(\mathbb{N})$, then for any $y > 0$,*

$$\mathbb{P}_x(f(X_t) - \mathbb{E}_x[f(X_t)] \geq y) \leq \exp\left(-\frac{y\mathrm{e}^{tK}}{2\|f\|_{\mathrm{Lip}_d}} \log\left(1 + \frac{yK}{\sinh(tK)\|\lambda + \nu\|_\infty \|f\|_{\mathrm{Lip}_d}}\right)\right).$$

*If $K = 0$, the latter inequality should be replaced by its limit as $K \to 0$.*

In Corollary 4.7 below, no particular boundedness assumption is made on the transition rates of the generator of the birth–death process $(X_t)_{t\geq 0}$.

**Corollary 4.7.** *Assume that the transition rates $\lambda$ and $\nu$ are, respectively, non-increasing and non-decreasing. If $f \in \mathrm{Lip}_d(\mathbb{N})$ and, further, $\|\Gamma f\|_\infty < +\infty$, then for any $y > 0$,*

$$\mathbb{P}_x(f(X_t) - \mathbb{E}_x[f(X_t)] \geq y) \leq \exp\left(-\frac{y}{2\|f\|_{\mathrm{Lip}_d}} \log\left(1 + \frac{y\|f\|_{\mathrm{Lip}_d}}{2t\|\Gamma f\|_\infty}\right)\right).$$

*Remark 4.8.* We are not able to extend to the stationary distribution $\pi$ the two previous deviation inequalities as $t$ goes to infinity. This is due to the non-positivity of the curvatures of the process $(X_t)_{t\geq 0}$ (cf. Remarks 4.2 and 4.5). In particular, it excludes the $M/M/\infty$ queueing process recently investigated by Chafaï [4] whose stationary distribution is the Poisson measure on $\mathbb{N}$. Therefore, we expect to recover the classical deviation inequality for the Poisson distribution by taking the limit as $t \to +\infty$ in an appropriate deviation estimate satisfied by the $M/M/\infty$ queueing process; an interesting problem of this nature will be addressed in a subsequent paper.

Our present purpose is to refine Corollaries 4.6 and 4.7 when the state space is finite, in order to establish, by a limiting argument, Poisson-type deviation estimates for the stationary distribution $\pi$. To do so, the crucial point is to obtain positive lower bounds on the curvatures in Propositions 4.1 and 4.4.

Our estimates below may be compared to that established by Houdré and Tetali [9], Proposition 4, under reversibility assumptions and without notions of curvature.

**Corollary 4.9.** *Assume that there exists $K > 0$ such that the transition rates $\lambda$ and $\nu$ satisfy $\min_{x \in \{1,\ldots,n\}} \lambda_{x-1} - \lambda_x + \nu_x - \nu_{x-1} \geq K$. If $f \in \mathrm{Lip}_d(\{0,1,\ldots,n\})$, then for any $y > 0$,*

$$\mathbb{P}_x(f(X_t) - \mathbb{E}_x[f(X_t)] \geq y) \leq \exp\left(-\frac{y}{2\|f\|_{\mathrm{Lip}_d}} \log\left(1 + \frac{2Ky}{(1-\mathrm{e}^{-2Kt})\|\lambda+\nu\|_\infty \|f\|_{\mathrm{Lip}_d}}\right)\right).$$

*In particular, letting $t \to +\infty$ above yields, under the stationary distribution $\pi$,*

$$\pi(f - \pi(f) \geq y) \leq \exp\left(-\frac{y}{2\|f\|_{\mathrm{Lip}_d}} \log\left(1 + \frac{2Ky}{\|\lambda+\nu\|_\infty \|f\|_{\mathrm{Lip}_d}}\right)\right).$$

Under different assumptions on the generator, we obtain the following, somewhat similar, estimate.

**Corollary 4.10.** *Assume that there exists $\rho > 0$ such that the transition rates $\lambda$ and $\nu$ satisfy $\min_{x \in \{1,\ldots,n-1\}} \min\{\lambda_{x-1} - \lambda_x, \nu_{x+1} - \nu_x\} \geq \rho$. If $f \in \mathrm{Lip}_d(\{0,1,\ldots,n\})$, then for any $y > 0$,*

$$\mathbb{P}_x(f(X_t) - \mathbb{E}_x[f(X_t)] \geq y) \leq \exp\left(-\frac{y}{2\|f\|_{\mathrm{Lip}_d}} \log\left(1 + \frac{2\rho y}{(1-\mathrm{e}^{-2\rho t})(\lambda_0 + \nu_n)\|f\|_{\mathrm{Lip}_d}}\right)\right).$$

*As $t \to +\infty$, we obtain the deviation inequality*

$$\pi(f - \pi(f) \geq y) \leq \exp\left(-\frac{y}{2\|f\|_{\mathrm{Lip}_d}} \log\left(1 + \frac{2\rho y}{(\lambda_0 + \nu_n)\|f\|_{\mathrm{Lip}_d}}\right)\right).$$

As an application of Corollary 4.9, let us recover the classical Gaussian deviation inequality for a Brownian-driven Ornstein–Uhlenbeck process constructed as a fluid limit of rescaled continuous-time Ehrenfest chains.

**Corollary 4.11.** *Let $(U_t)_{t \geq 0}$ be the Brownian-driven Ornstein–Uhlenbeck process given by*

$$U_t = z_0 \mathrm{e}^{-t} + \sqrt{2\lambda\nu} \int_0^t \mathrm{e}^{-(t-s)} \mathrm{d}B_s, \qquad t > 0,$$

*where $z_0 \in \mathbb{R}$ and the positive parameters $\lambda$ and $\nu$ are such that $\lambda + \nu = 1$. Then for any Lipschitz function $f$ on $\mathbb{R}$ with Lipschitz constant $\|f\|_{\mathrm{Lip}}$, the classical Gaussian deviation inequality holds:*

$$\mathbb{P}_{z_0}(f(U_t) - \mathbb{E}_{z_0}[f(U_t)] \geq y) \leq \exp\left(-\frac{y^2}{(1-\mathrm{e}^{-2t})\nu \|f\|_{\mathrm{Lip}}^2}\right), \qquad y > 0.$$

**Proof.** Let $(X_t^n)_{t\geq 0}$ be the continuous-time Ehrenfest chain on $\{0,1,\ldots,n\}$ starting from some $x_n \in \{0,1,\ldots,n\}$ and with generator given by

$$\mathcal{L}_n f(x) = \lambda(n-x)(f(x+1) - f(x)) + \nu x(f(x-1) - f(x)), \qquad x \in \{0,1,\ldots,n\}.$$

Suppose that $\lim_{n\to+\infty} x_n/n = \lambda$ and define the process $(Z_t^n)_{t\geq 0}$ by $Z_t^n = (X_t^n - \lambda n)/\sqrt{n}$, $t > 0$. Further assume that the sequence of initial states $(Z_0^n)_{n\in\mathbb{N}}$ converges to $z_0$. By the central limit theorem in Ethier and Kurtz [7], Chapter 11, the sequence $(Z_t^n)_{t\geq 0}$ converges to the process $(U_t)_{t\geq 0}$ as $n$ goes to infinity.

Now, fix $n \in \mathbb{N}\setminus\{0\}$, $t > 0$, and consider the function $h_n = f \circ \phi_n$, where $\phi_n$ is defined on $\{0,1,\ldots,n\}$ by $\phi_n(x) = (x - n\lambda)/\sqrt{n}$. Then $h_n \in \mathrm{Lip}_d(\{0,1,\ldots,n\})$ with Lipschitz constant at most $n^{-1/2}\|f\|_{\mathrm{Lip}}$. Therefore, we can apply Corollary 4.9 to $(X_t^n)_{t\geq 0}$ and $h_n$, with $K = \lambda + \nu = 1$, to obtain, for any $y > 0$,

$$\mathbb{P}_{x_n}(h_n(X_t^n) - \mathbb{E}_{x_n}[h_n(X_t^n)] \geq y) \leq \exp\left(-\frac{y\sqrt{n}}{2\|f\|_{\mathrm{Lip}}} \log\left(1 + \frac{2y}{(1-\mathrm{e}^{-2t})\sqrt{n}\nu\|f\|_{\mathrm{Lip}}}\right)\right).$$

Finally, letting $n$ tend to infinity in the above inequality yields the result. $\square$

## 4.3. A multidimensional deviation inequality for the $M/M/1$ queue

In this section, we give a Poisson-type deviation estimate for a sample of the $M/M/1$ queueing process. It is an irreducible birth–death process $(X_t)_{t\geq 0}$ whose generator is given by

$$\mathcal{L}f(x) = \lambda(f(x+1) - f(x)) + \nu 1_{\{x\neq 0\}}(f(x-1) - f(x)), \qquad x \in \mathbb{N},$$

where the numbers $\lambda$ and $\nu$ are positive. The existence of an integration by parts formula for the associated semigroup, together with a tensorization procedure of the Laplace transform, allow us to provide, in Corollary 4.12 below, a multidimensional deviation inequality for the $M/M/1$ queue.

We say in the sequel that a function $f:\mathbb{N}^n \to \mathbb{R}$ is $\ell^1$-*Lipschitz* if

$$\|f\|_{\mathrm{Lip}(n)} = \sup_{x\neq y} \frac{|f(x) - f(y)|}{\|x-y\|_1} < +\infty,$$

where $\|\cdot\|_1$ denotes the $\ell^1$-norm $\|z\|_1 = \sum_{i=1}^n |z_i|$, $z \in \mathbb{N}^n$.

**Corollary 4.12.** *Consider the sample $X^n = (X_{t_1},\ldots,X_{t_n})$, $0 = t_0 < t_1 < \cdots < t_n = T$, and let $f$ be $\ell^1$-Lipschitz on $\mathbb{N}^n$. Then, for any $y > 0$,*

$$\mathbb{P}_x(f(X^n) - \mathbb{E}_x[f(X^n)] \geq y) \leq \exp\left(-\frac{y}{2n\|f\|_{\mathrm{Lip}(n)}} \log\left(1 + \frac{y}{Tn(\lambda+\nu)\|f\|_{\mathrm{Lip}(n)}}\right)\right). \tag{4.5}$$

**Proof.** Let $u$ be a Lipschitz function on $\mathbb{N}$ with Lipschitz constant $\|u\|_{\text{Lip}(1)}$ and let $t > 0$. Rewriting the proof of Theorem 3.4 for the $M/M/1$ queue yields, for any $\tau > 0$,

$$\mathbb{E}_x[\mathrm{e}^{\tau u(X_t)}] \leq \exp\{\tau \mathbb{E}_x[u(X_t)] + h(\tau, t, \|u\|_{\text{Lip}(1)})\}, \tag{4.6}$$

where $h$ is the function defined on $(\mathbb{R}_+)^3$ by $h(\tau, t, z) = t(\lambda + \nu)(\mathrm{e}^{\tau z} - \tau z - 1)$.

To obtain a multidimensional version of (4.6), the idea is to tensorize the Laplace transform with respect to the $\ell^1$-metric via an integration by parts formula satisfied by the semigroup $(P_t)_{t \geq 0}$ of the $M/M/1$ queueing process.

First, observe that we have the commutation relation $\mathcal{L}\mathrm{d}^+u = \mathrm{d}^+\mathcal{L}u$, where $\mathrm{d}^+$ is the forward gradient $\mathrm{d}^+u(x) = u(x+1) - u(x)$, $x \in \mathbb{N}$. It implies $P_t \mathrm{d}^+ u = \mathrm{d}^+ P_t u$ for any $t \geq 0$, which, in turn, yields the integration by parts formula:

$$\sum_{y \in \mathbb{N}} u(y) P_t(x+1, y) = \sum_{y \in \mathbb{N}} u(y+1) P_t(x, y), \qquad x \in \mathbb{N}. \tag{4.7}$$

Let $f$ be $\ell^1$-Lipschitz on $\mathbb{N}^n$. Set $f_n := f$ and define, for any $k = 1, \ldots, n-1$, the function $f_k$ on $\mathbb{N}^k$ by

$$f_k(x_1, \ldots, x_k)$$
$$:= \sum_{x_{k+1}, \ldots, x_n \in \mathbb{N}} f(x_1, \ldots, x_k, x_{k+1}, \ldots, x_n) P_{t_{k+1} - t_k}(x_k, x_{k+1}) \cdots P_{t_n - t_{n-1}}(x_{n-1}, x_n).$$

Let $x_1, \ldots, x_{k-1}, y \in \mathbb{N}$. Recursively using (4.7), we have

$$f_k(x_1, \ldots, x_{k-1}, y+1)$$
$$= \sum_{x_{k+1}, \ldots, x_n \in \mathbb{N}} f(x_1, \ldots, x_{k-1}, y+1, x_{k+1}, x_{k+2}, \ldots, x_n) P_{t_{k+1} - t_k}(y+1, x_{k+1})$$
$$\times P_{t_{k+2} - t_{k+1}}(x_{k+1}, x_{k+2}) \cdots P_{t_n - t_{n-1}}(x_{n-1}, x_n)$$
$$= \sum_{x_{k+1}, \ldots, x_n \in \mathbb{N}} f(x_1, \ldots, x_{k-1}, y+1, x_{k+1}+1, x_{k+2}, \ldots, x_n) P_{t_{k+1} - t_k}(y, x_{k+1})$$
$$\times P_{t_{k+2} - t_{k+1}}(x_{k+1}+1, x_{k+2}) \cdots P_{t_n - t_{n-1}}(x_{n-1}, x_n)$$
$$= \cdots$$
$$= \sum_{x_{k+1}, \ldots, x_n \in \mathbb{N}} f(x_1, \ldots, x_{k-1}, y+1, x_{k+1}+1, \ldots, x_{n-1}+1, x_n)$$
$$\times P_{t_{k+1} - t_k}(y, x_{k+1}) \cdots P_{t_{n-1} - t_{n-2}}(x_{n-2}, x_{n-1}) P_{t_n - t_{n-1}}(x_{n-1}+1, x_n)$$
$$= \sum_{x_{k+1}, \ldots, x_n \in \mathbb{N}} f(x_1, \ldots, x_{k-1}, y+1, x_{k+1}+1, \ldots, x_{n-1}+1, x_n+1)$$
$$\times P_{t_{k+1} - t_k}(y, x_{k+1}) \cdots P_{t_{n-1} - t_{n-2}}(x_{n-2}, x_{n-1}) P_{t_n - t_{n-1}}(x_{n-1}, x_n).$$

Hence, we obtain, for any $k = 1, \ldots, n$ and any $x_1, \ldots, x_{k-1} \in \mathbb{N}$,

$$\|f_k(x_1, \ldots, x_{k-1}, \cdot)\|_{\mathrm{Lip}(1)}$$
$$= \sup_{y \in \mathbb{N}} |f_k(x_1, \ldots, x_{k-1}, y+1) - f_k(x_1, \ldots, x_{k-1}, y)|$$
$$\leq \sup_{y \in \mathbb{N}} \sum_{x_{k+1}, \ldots, x_n \in \mathbb{N}} |f(x_1, \ldots, x_{k-1}, y+1, \ldots, x_n+1) - f(x_1, \ldots, x_{k-1}, y, \ldots, x_n)|$$
$$\times P_{t_{k+1}-t_k}(y, x_{k+1}) \cdots P_{t_{n-1}-t_{n-2}}(x_{n-2}, x_{n-1}) P_{t_n-t_{n-1}}(x_{n-1}, x_n)$$
$$\leq (n-k+1)\|f\|_{\mathrm{Lip}(n)}$$
$$\leq n \|f\|_{\mathrm{Lip}(n)}. \tag{4.8}$$

Successively using the inequality (4.6) in the following lines with the one-dimensional Lipschitz functions $x_k \mapsto f_k(*, x_k)$, $k = n, n-1, \ldots, 1$, and substituting the upper bound of (4.8) into the right-hand side of (4.6) (since the function $h$ is non-decreasing in its last variable), we have

$$\mathbb{E}_x[e^{\tau f(X^n)}]$$
$$= \sum_{x_1, \ldots, x_{n-1} \in \mathbb{N}} \sum_{x_n \in \mathbb{N}} \exp\{\tau f_n(x_1, \ldots, x_n)\} P_{t_n-t_{n-1}}(x_{n-1}, x_n) \cdots P_{t_1}(x, x_1)$$
$$\leq \exp\{h(\tau, t_n - t_{n-1}, n\|f\|_{\mathrm{Lip}(n)})\}$$
$$\times \sum_{x_1, \ldots, x_{n-2} \in \mathbb{N}} \sum_{x_{n-1} \in \mathbb{N}} \exp\{\tau f_{n-1}(x_1, \ldots, x_{n-1})\} P_{t_{n-1}-t_{n-2}}(x_{n-2}, x_{n-1}) \cdots P_{t_1}(x, x_1)$$
$$\leq \exp\{h(\tau, t_n - t_{n-1}, n\|f\|_{\mathrm{Lip}(n)}) + h(\tau, t_{n-1} - t_{n-2}, n\|f\|_{\mathrm{Lip}(n)})\}$$
$$\times \sum_{x_1, \ldots, x_{n-3} \in \mathbb{N}} \sum_{x_{n-2} \in \mathbb{N}} \exp\{\tau f_{n-2}(x_1, \ldots, x_{n-2})\} P_{t_{n-2}-t_{n-3}}(x_{n-3}, x_{n-2}) \cdots P_{t_1}(x, x_1)$$
$$\leq \cdots$$
$$\leq \exp\left(\sum_{k=1}^{n-1} h(\tau, t_{n-k+1} - t_{n-k}, n\|f\|_{\mathrm{Lip}(n)})\right) \sum_{x_1 \in \mathbb{N}} e^{\tau f_1(x_1)} P_{t_1}(x, x_1)$$
$$\leq \exp\left(\sum_{k=1}^{n} h(\tau, t_{n-k+1} - t_{n-k}, n\|f\|_{\mathrm{Lip}(n)})\right) \exp\left\{\tau \sum_{x_1 \in \mathbb{N}} f_1(x_1) P_{t_1}(x, x_1)\right\}$$
$$= \exp\left(\sum_{k=1}^{n} h(\tau, t_k - t_{k-1}, n\|f\|_{\mathrm{Lip}(n)})\right) \exp\{\tau \mathbb{E}_x[f(X^n)]\}$$
$$= \exp\{\tau \mathbb{E}_x[f(X^n)] + T(\lambda + \nu)(\exp\{\tau n\|f\|_{\mathrm{Lip}(n)}\} - \tau n\|f\|_{\mathrm{Lip}(n)} - 1)\}.$$

Dividing by $e^{\tau \mathbb{E}_x[f(X^n)]}$ and using Chebyshev's inequality completes the proof. $\square$

***Remark 4.13.*** To conclude this work, note that Corollary 4.12 does not allow us to extend the deviation inequality (4.5) to functionals on path spaces. Thus it would be interesting to suitably refine such an estimate in terms of the increments $\Delta_i = t_i - t_{i-1}$, as $\Delta_i \to 0$.